\documentclass[12pt]{amsart}%
\usepackage{amsfonts}
\usepackage{amsmath}
\usepackage{amssymb}
\usepackage{amsthm}
\usepackage{booktabs}
\usepackage{mathtools}
\usepackage{commath}
\usepackage{cleveref}
\usepackage{enumerate}
\usepackage{bm}

\makeatletter
\newtheorem*{rep@theorem}{\rep@title}
\newcommand{\newreptheorem}[2]{%
\newenvironment{rep#1}[1]{%
 \def\rep@title{#2 \ref{##1}}%
 \begin{rep@theorem}}%
 {\end{rep@theorem}}}
\makeatother

\newtheorem{theorem}{Theorem}[section]
\newreptheorem{theorem}{Theorem}

\newtheorem{corollary}[theorem]{Corollary}
\newtheorem{lemma}[theorem]{Lemma}
\newreptheorem{corollary}{Corollary}

\newtheorem{proposition}[theorem]{Proposition}

\theoremstyle{definition}
\newtheorem{definition}[theorem]{Definition}

\newtheorem{remark}[theorem]{Remark}
\newtheorem{example}[theorem]{Example}

\newcommand{\bC}{\mathbb{C}}
\newcommand{\bD}{\mathbb{D}}
\newcommand{\bQ}{\mathbb{Q}}
\newcommand{\bR}{\mathbb{R}}
\newcommand{\bZ}{\mathbb{Z}}
\newcommand{\cA}{\mathcal{A}}

\newcommand{\cD}{\mathcal{D}}

\newcommand{\cH}{\mathcal{H}}
\newcommand{\cI}{\mathcal{I}}
\newcommand{\cJ}{\mathcal{J}}

\newcommand{\cSD}{\mathcal{SD}}

\newcommand{\fB}{\mathfrak{B}}
\newcommand{\fD}{\mathfrak{D}}
\newcommand{\fF}{\mathfrak{F}}

\newcommand{\fS}{\mathfrak{S}}
\newcommand{\fZ}{\mathfrak{Z}}

\DeclareMathOperator{\bideg}{bideg}
\DeclareMathOperator{\End}{End}
\DeclareMathOperator{\Cliff}{Cliff}
\DeclareMathOperator{\GL}{GL}
\DeclareMathOperator{\GrFrob}{GrFrob}
\DeclareMathOperator{\Hilb}{Hilb}
\DeclareMathOperator{\Hom}{Hom}
\DeclareMathOperator{\id}{id}
\DeclareMathOperator{\perm}{perm}
\DeclareMathOperator{\rS}{S}
\DeclareMathOperator{\sgn}{sgn}
\DeclareMathOperator{\Span}{Span}
\DeclareMathOperator{\Sym}{Sym}
\DeclareMathOperator{\Tr}{Tr}

\newcommand{\res}{\!\!\downarrow}

\title[Harmonic differential forms I. Semi-invariants]{Harmonic differential forms for pseudo-reflection groups \\ I. Semi-invariants}
\author{Joshua P. Swanson and Nolan R. Wallach}
\date{\today}

\begin{document}
  
\begin{abstract}
  We give a type-independent construction of an explicit basis for the
  semi-invariant harmonic differential forms of an arbitrary pseudo-reflection group
  in characteristic zero.
  Our ``top-down'' approach uses the methods of Cartan's exterior calculus and is
  in some sense dual to related work of Solomon \cite{MR0154929},
  Orlik--Solomon \cite{MR575083},
  and Shepler \cite{MR1714136,MR2163712} describing (semi-)invariant differential forms.
  We apply our results to a recent conjecture of
  Zabrocki \cite{1902.08966} which provides a representation theoretic-model
  for the Delta conjecture of Haglund--Remmel--Wilson \cite{MR3811519} in
  terms of a certain non-commutative coinvariant algebra for the symmetric
  group. In particular, we verify the alternating component of a
  specialization of Zabrocki's conjecture.
\end{abstract}
\maketitle

\section{Introduction}

Recently there has been a great deal of research activity in algebraic
combinatorics studying diagonal actions of the symmetric
group $\fS_n$ on $k$ sets of $n$ commuting indeterminants and $\ell$ sets of
$n$ anti-commuting indeterminants. Orellana--Zabrocki \cite{1906.01125}
describe the $\fS_n$-invariants of these polynomial rings combinatorially
and summarize some of their history. The $k=2, \ell=0$ case has received a very large
amount of attention through the study of the \textit{diagonal coinvariants}
$\bQ[\bm{x}_n, \bm{y}_n]/\cD_n$ where $\cD_n$ is the \textit{diagonal coinvariant ideal}
generated by all homogeneous $\fS_n$-invariants of positive degree and $\bm{x}_n$
is shorthand for $x_1, \ldots, x_n$ \cite{MR1918676,MR1256101}. Zabrocki
\cite{1902.08966} has recently given a conjectured description of the
tri-graded $\fS_n$-module isomorphism type of the coinvariant algebra
when $m=2, \ell=1$:
\begin{equation}\label{eq:zab}
  \GrFrob(\bQ[\bm{x}_n, \bm{y}_n, \bm{\theta}_n]/\cSD_n; q, t, z)
      = \sum_{i=1}^n z^{n-k} \Delta_{e_{k-1}}' e_n,
\end{equation}
where $\GrFrob$ is the tri-graded \textit{Frobenius series},
$\cSD_n$ is the super diagonal coinvariant ideal generated by
homogeneous $\fS_n$-invariants of positive degree,
$e_n$ is an elementary symmetric function, and $\Delta_f'$ is a
certain modified Macdonald eigenoperator. \Cref{eq:zab} may be interpreted as a conjectural
representation-theoretic model for the Delta conjecture of Haglund--Remmel--Wilson
\cite{MR3811519}, and both conjectures remain open. 
See \cite{MR3811519,1902.08966} for details and further references. 

The special case $t=0$ of these conjectures involving one set of commuting and
one set of anti-commuting variables has received special attention.
Haglund--Rhoades--Shimozono
\cite{MR3783430} had earlier given a different representation-theoretic model
for this specialization of the Delta conjecture. Motivated by Zabrocki's conjecture,
the second author \cite{1906.11787} gave a conjectural description of the
\textit{harmonics} of $\bQ[\bm{x}_n, \bm{\theta}_n]/\cJ_n$ where
$\cJ_n$ is the super coinvariant ideal generated by homogeneous
$\fS_n$-invariants of positive degree.
This description further motivated Rhoades--Wilson
\cite{1906.03315} to recently construct another representation-theoretic model for the $t=0$
specialization of the Delta conjecture arising from the leading terms of
those harmonics. In this paper, we restrict our attention to
natural pseudo-reflection group generalizations of the $t=0$ case of
Zabrocki's coinvariant algebra model.

In \cite{1906.11787}, several of the implications of Zabrocki's conjecture at $t=0$
were proven, including a determination of the bi-graded Hilbert series
of the alternants in the $\fS_n$-coinvariants and certain vanishing bounds.
The key to the results in \cite{1906.11787} is the fact that one is looking
at differential forms with polynomial coefficients and
thus one has at one's disposal the full power of Cartan's exterior calculus.
The idea of adding anti-commuting variables to prove theorems about the
commuting variables appeared in the work of Solomon \cite{MR0154929} and,
later, in the paper of Orlik--Solomon \cite{MR575083} in the case of a finite unitary subgroup
$G \leq \GL_n(\bC)$ generated by pseudo-reflections. Differential forms and derivations
have also been exploited in this context more recently \cite{MR3959046,RSS19}.
The purpose of the present paper is to use the methods of differential forms to
give a uniform generalization of the type $A$ description in \cite{1906.11787} of the
alternant polynomial differential forms and coinvariants for an
arbitrary pseudo-reflection group.

Our main results are as follows.
See the subsequent sections for missing definitions and
\Cref{sec:super_alternants} and \Cref{sec:coinvariant_alternants} for proofs.
Let $V$ be an $n$-dimensional vector space over an arbitrary field $F$ of characteristic $0$,
let $G \leq \GL(V)$ be a pseudo-reflection group, let $M$ be an $r$-dimensional
$G$-module, and let $\chi$ be a one-dimensional character of $G$.

We consider the
\textit{semi-invariant differential forms}
\begin{equation}\label{eq:semis}
  (\rS(V^*) \otimes \wedge M^*)^\chi
\end{equation}
where $\rS(V^*)$ is the symmetric algebra on $V^*$, $\wedge M^*$ is the
exterior algebra on $M^*$, and $W^\chi \coloneqq \{w \in W : \sigma \cdot w
= \chi(\sigma)w, \forall \sigma \in G\}$.
In certain circumstances, such as when $M=V$, we give two explicit
bases for \eqref{eq:semis}. The following result may be thought
of as an analogue of Solomon's classic result \cite{MR0154929} describing
$(\rS(V^*) \otimes \wedge V^*)^G$ as a Grassmann algebra over $\rS(V^*)^G$. Here
$\Delta_\chi$ is the minimal-degree
element in the $\chi$-isotypic component of $\rS(V^*)$, the $d_i^*$ are certain
differential operators depending on $M$ (see \Cref{def:d_i_star}), the $f_i$ are
basic invariants of $G$, and $J_{M^*}$ is the \textit{Jacobian} of $M^*$
(see \Cref{def:jacobian}).

\begin{reptheorem}{thm:super_alternant_basis.2}
  Suppose $J_{M^*} \mid \Delta_\chi$. Then either of the sets
    \[ \{f_1^{a_1} \cdots f_n^{a_n} d_{i_1}^* \cdots d_{i_k}^* \Delta_\chi
           : 1 \leq i_1 < \cdots < i_k \leq r, a_j \in \bZ_{\geq 0}\} \]
  or
    \[ \{d_{i_1}^* \cdots d_{i_k}^* f_1^{a_1} \cdots f_n^{a_n} \Delta_\chi
           : 1 \leq i_1 < \cdots < i_k \leq r, a_j \in \bZ_{\geq 0}\} \]
  form bases for $(\rS(V^*) \otimes \wedge M^*)^\chi$.
\end{reptheorem}

We have the following enumerative corollary. Let $e_i^{M^*}$ denote the
\textit{exponents} of $M^*$ (see \Cref{def:exponents}) and let $d_i$ denote
the \textit{degrees} of $G$.

\begin{repcorollary}{cor:super_alternant_basis.2}
  Suppose $J_{M^*} \mid \Delta_\chi$. Then
    \[ \Hilb((\rS(V^*) \otimes \wedge M^*)^\chi; q, z)
        = q^{\deg \Delta_\chi} \frac{\prod_{i=1}^r (1 + zq^{-e_i^{M^*})}}
           {\prod_{i=1}^n (1 - q^{d_i})}. \]
\end{repcorollary}

\begin{remark}
  Orlik--Solomon \cite{MR575083} described $(\rS(V^*) \otimes \wedge M^*)^G$
  when $J_M = \Delta_M$, up to a non-zero scalar, as an exterior algebra over
  $\rS(V^*)^G$, and Shepler \cite{MR1714136,MR2163712} gave an analogous
  result for $(\rS(V^*) \otimes \wedge V^*)^\chi$ using ``$\chi$-wedging.''
  These exterior algebra structures may be considered ``bottom-up''
  descriptions. \Cref{thm:super_alternant_basis.2} involves the operators $d_i^*$,
  which decrease $\rS(V^*)$-degree, applied to $\Delta_\chi$, so this and related work
  in \cite[\S6]{MR2163712} can be considered ``top-down.'' The top-down description turns
  out to be more useful in our analysis of the coinvariant algebra.
\end{remark}

Let $\cJ_M^*$ denote the ideal in $\rS(V^*) \otimes \wedge M^*$ generated by all
homogeneous $G$-invariants of positive degree. Our overarching motivation has been
to describe the semi-invariant elements of the coinvariant algebra,
\begin{equation}\label{eq:semi_coinvs}
  (\rS(V^*) \otimes \wedge M^*/\cJ_M^*)^\chi.
\end{equation}
We give the following explicit basis for \eqref{eq:semi_coinvs} using the
\textit{harmonics} of $\rS(V^*) \otimes \wedge M^*$ (see \Cref{def:harmonics}).

\begin{reptheorem}{thm:super_coinvariant_basis}
  Suppose $J_{M^*} \mid \Delta_\chi$ and $M^G = 0$. Then the $2^r$ elements
    \[ \{d_{i_1}^* \cdots d_{i_k}^* \Delta_\chi \mid 1 \leq i_1 < \cdots < i_k \leq r\} \]
  form a basis of $\cH(\rS(V^*) \otimes \wedge M^*)^\chi$,
  and their images descend to a basis of
  $(\rS(V^*) \otimes \wedge M^*/\cJ_M^*)^\chi$.
\end{reptheorem}

\begin{repcorollary}{cor:super_coinvariant_basis}
  Suppose $J_{M^*} \mid \Delta_\chi$ and $M^G = 0$. Then
  \begin{equation}\label{eq:cor_super_coinvariants}
    \Hilb((\rS(V^*) \otimes \wedge M^*/\cJ_M^*)^\chi; q, z)
        = q^{\deg \Delta_\chi}\prod_{i=1}^r (1 + zq^{-e_i^{M^*}}).
  \end{equation}
\end{repcorollary}

The hypotheses of \Cref{thm:super_coinvariant_basis} are satisfied whenever
$\chi = \det_{M^*}$, $M^G = 0$, and the pseudo-reflections of $G$ act by pseudo-reflections
or the identity on $M$. In this case, the right-hand side of \eqref{eq:cor_super_coinvariants}
is
  \[ \prod_{i=1}^r (q^{e_i^{M^*}} + z). \]
The alternating component of the
$t=0$ specialization of Zabrocki's conjecture then follows from
\Cref{cor:super_coinvariant_basis} when $G$ consists of $n \times n$ permutation matrices
and $M$ is the $(n-1)$-dimensional standard representation.

The classical coinvariant algebra of $G$ is $\rS(V^*)/\cI^*$ where $\cI^*$ is
the ideal generated by all non-constant homogeneous $G$-invariants. It is
well-known that the top-degree component of $\rS(V^*)/\cI^*$ is the image of
$\rS(V^*)^{\det_V}$, which has motivated much of our work.
We will explore which bidegrees of
$\rS(V^*) \otimes \wedge M^*/\cJ_M^*$ are non-zero in a future article
\cite{sw19.2}.

The rest of the paper is organized as follows.
In \Cref{sec:pairings}, we review background material on differential forms and
pseudo-reflection groups. \Cref{sec:exps} concerns exponents and basic derivations
of pseudo-reflection groups. \Cref{sec:super_alternants} introduces
$G$-harmonics for the polynomial algebra and constructs explicit bases for
certain pieces of the differential algebras; see
\Cref{thm:super_alternant_basis} and \Cref{thm:super_alternant_basis.2}.
In \Cref{sec:coinvariant_alternants}, we prove
our main result, \Cref{thm:super_coinvariant_basis}, which gives a basis for the harmonics
and the coinvariants. In \Cref{sec:Jm_Dm}, we discuss the technical
condition $\Delta_M = J_M$, up to a non-zero scalar, which appears in some of our results.
A more explicit but less general version of many of these results
that may be more palatable to algebraic combinatorialists can be found in
\cite{1908.00196}.

\section{Polynomial differential forms}\label{sec:pairings}

We now describe several actions and pairings involving
polynomial differential forms and related objects. All of the constructions
in this section use standard ideas from differential geometry.

Let $F$ be a field of characteristic $0$ and let $V$ be an $F$-vector space of
dimension $n < \infty$. We identify $V^{**} = V$. Let $G$ be an arbitrary subgroup of
$\GL(V)$ and let $M$ be an $F[G]$-module of dimension $r < \infty$.

\subsection{Symmetric and exterior algebras}

\begin{definition}
  Let $\rS(V^*) \coloneqq \Sym(V^*)$ denote the symmetric algebra on $V^*$ over
  $F$, i.e.~the algebra of polynomial functions on $V$. Let
  $\wedge M^*$ denote the exterior or \textit{Grassmann} algebra on $M^*$.
\end{definition}

The tensor products $\rS(V) \otimes \wedge M$ and $\rS(V^*) \otimes \wedge M^*$ are
algebras of differential forms with polynomial coefficients.
Each of $\rS(V)$, $\rS(V^*)$, $\wedge M$, and $\wedge M^*$ is naturally graded,
so we have four non-commutative, bigraded $F$-algebras
\begin{equation}\label{eq:four_algebras}
  \rS(V) \otimes \wedge M, \qquad
  \rS(V) \otimes \wedge M^*, \qquad
  \rS(V^*) \otimes \wedge M, \qquad
  \rS(V^*) \otimes \wedge M^*.
\end{equation}

The $G$-action on $V$ extends multiplicatively to yield natural $G$-actions on
$\rS(V)$ and $\wedge M$. As usual, $G$ acts on $V^*$
\textit{contragrediently} via $g^*\lambda \coloneqq \lambda \circ g^{-1}$.
Thus $\rS(V)$, $\rS(V^*)$, $\wedge M$, and $\wedge M^*$ are all naturally
graded $G$-modules, so the algebras in \eqref{eq:four_algebras} are
bigraded $G$-modules via the diagonal actions of
  \[ g \otimes g, \qquad
      g \otimes g^*, \qquad
      g^* \otimes g, \qquad
      g^* \otimes g^*. \]

\subsection{Pairings and differential operators on $\rS(V^*) \otimes \wedge M^*$}\label{ssec:pairings}

Our next goal is to describe natural actions of $\rS(V) \otimes \wedge M$ and
$\rS(V) \otimes \wedge M^*$ on $\rS(V^*) \otimes \wedge M^*$. These
actions will be fundamental in later sections.

The following observation is essentially trivial.

\begin{lemma}
  The natural pairing $M ^*\times M \to F$ given by $\langle \lambda, v\rangle
  \coloneqq \lambda(v)$ is $G$-invariant and perfect.
\end{lemma}

By ``perfect'', we mean that $\langle \lambda, v\rangle = 0$ for all $v \in V$ implies
$\lambda = 0$, and $\langle \lambda, v\rangle = 0$ for all $\lambda \in V^*$ implies
$v=0$. We may naturally extend such pairings to symmetric, exterior, and tensor products as
follows.

\begin{lemma}\label{lem:amp_pairings.1}
  Suppose $W, W_1, W_2$ are finite-dimensional $G$-modules with
  $G$-invariant perfect pairings
    \[ \langle -, -\rangle \colon W^* \times W \to F \]
  and
    \[ \langle -, -\rangle_i \colon W_i^* \times W_i \to F. \]
  \begin{enumerate}[(i)]
    \item The pairing
      \begin{align*}
        W_1^* \otimes W_2^* \times W_1 \otimes W_2 &\to F \\
        \langle \omega_1 \otimes \omega_2, w_1 \otimes w_2\rangle
          &\coloneqq \langle \omega_1, w_1\rangle_1 \langle \omega_2, w_2\rangle_2
      \end{align*}
      is $G$-invariant and perfect.
    \item The pairing
      \begin{align*}
        \wedge W^* \times \wedge W &\to F \\
        \langle \omega_1 \wedge \cdots \wedge \omega_k,
            w_1 \wedge \cdots \wedge w_\ell\rangle
          &\coloneqq \delta_{k\ell} \det(\langle \omega_i, w_j\rangle)_{i,j=1}^k \\
          &\coloneqq \delta_{k\ell} \sum_{\sigma \in \fS_k} (-1)^{\sgn(\sigma)}
               \prod_{i=1}^k \langle \omega_i, w_{\sigma(i)}\rangle
      \end{align*}
      is $G$-invariant and perfect.
    \item The pairing
      \begin{align*}
        \rS(W^*) \times \rS(W) &\to F \\
        \langle \omega_1 \cdots \omega_k, w_1 \cdots w_\ell\rangle
          &\coloneqq \delta_{k\ell} \perm(\langle \omega_i, w_j\rangle)_{i,j=1}^k \\
          &\coloneqq \delta_{k\ell} \sum_{\sigma \in \fS_k}
               \prod_{i=1}^k \langle \omega_i, w_{\sigma(i)}\rangle
      \end{align*}
      is $G$-invariant and perfect, where $\perm$ denotes the matrix permanent.
  \end{enumerate}
  
  \begin{proof}
    In each case, $G$-invariance is immediate and non-degeneracy can be
    checked quickly using dual bases.
  \end{proof}
\end{lemma}

\begin{corollary}\label{cor:pairing}
  We have a natural $G$-invariant perfect pairing
  \begin{equation}\label{eq:perfect_pairing}
    \langle -, -\rangle \colon \rS(V^*) \otimes \wedge M^* \times
      \rS(V) \otimes \wedge M \to F.
  \end{equation}
  Furthermore, we have natural $G$-equivariant identifications of
  $\rS(V^*)$ with $\rS(V)^*$ and of $\wedge M^*$ with $(\wedge M)^*$.
\end{corollary}

\begin{definition}
  For $s \in \rS(V)$ and $f \in \rS(V^*)$, we have multiplication operators
  defined by
  \begin{alignat*}{4}
    m_s \colon \rS(V) &\to \rS(V) \qquad&&\text{and}\qquad&
    m_f \colon \rS(V^*) &\to \rS(V^*) \\
    m_s(u) &\coloneqq su \qquad&&\text{and}\qquad&
    m_f(h) &\coloneqq fh.
  \end{alignat*}
  For $m \in \wedge M$ and $\mu \in \wedge M^*$, we have
  multiplication operators
  \begin{alignat*}{4}
    \epsilon_m \colon {\wedge M} &\to \wedge M \qquad&&\text{and}\qquad&
    \epsilon_\mu \colon {\wedge M^*} &\to \wedge M^* \\
    \epsilon_m(\ell) &\coloneqq m \wedge \ell \qquad&&\text{and}\qquad&
    \epsilon_\mu(\nu) &\coloneqq \mu \wedge \nu.
  \end{alignat*}
\end{definition}

\begin{definition}
  For $s \in \rS(V)$ and $f \in \rS(V^*)$, we have adjoint operators defined by
  \begin{alignat*}{4}
    \partial_s \colon \rS(V^*) &\to \rS(V^*) \qquad&&\text{and}\qquad&
    \partial_f \colon \rS(V) &\to \rS(V) \\
    \langle \partial_s(h), u\rangle &= \langle h, m_s(u)\rangle
      \qquad&&\text{and}\qquad&
    \langle h, \partial_f(u)\rangle &= \langle m_f(h), u\rangle
  \end{alignat*}
  for all $u \in \rS(V)$, $h \in \rS(V^*)$.

  For $m \in \wedge M$ and $\mu \in \wedge M^*$, we have adjoint operators defined by
  \begin{alignat*}{4}
    \iota_m \colon {\wedge} M^* &\to \wedge M^* \qquad&&\text{and}\qquad&
    \iota_\mu \colon {\wedge} M &\to \wedge M \\
    \langle \iota_m(\nu), \ell\rangle &= \langle \nu, \epsilon_m(\ell)\rangle
      \qquad&&\text{and}\qquad&
    \langle \nu, \iota_\mu(\ell)\rangle &= \langle \epsilon_\mu(\nu), \ell\rangle
  \end{alignat*}
  for all $\ell \in \wedge M$, $\nu \in \wedge M^*$.
\end{definition}

Combining these operators yields actions of $\rS(V) \otimes \wedge M$
and $\rS(V) \otimes \wedge M^*$ on $\rS(V^*) \otimes \wedge M^*$ as follows.
These actions will be used prominently in later sections.

\begin{lemma}\label{lem:actions}
  The maps
  \begin{align*}
    \rS(V) \otimes \wedge M &\to \End_F(\rS(V^*) \otimes \wedge M^*) \\
    s \otimes m &\mapsto \partial_s \otimes \iota_m
  \end{align*}
  and
  \begin{align*}
      \rS(V) \otimes \wedge M^* &\to \End_F(\rS(V^*) \otimes \wedge M^*) \\
      s \otimes \mu &\mapsto \partial_s \otimes \epsilon_\mu
  \end{align*}
  are $G$-equivariant $F$-algebra morphisms.

  
  \begin{proof}
    The multiplication operators taken together yield $F$-algebra morphisms, and the
    same is true of their adjoints. By definition, $g \in G$ acts on
    $\psi \in \End_F(\rS(V^*) \otimes \wedge M^*)$ by
    $(g\psi)(\omega) \coloneqq g\psi(g^{-1}\omega)$. For $G$-equivariance,
    we see immediately that
    \begin{alignat*}{4}
      m_{gs}(gu) &= gm_s(u) \qquad&&\text{and}\qquad&
        m_{gf}(gh) &= gm_f(h) \\
      \epsilon_{gm}(gu) &= g\epsilon_m(u)
        \qquad&&\text{and}\qquad&
       \epsilon_{g\mu}(g\nu) &= g\epsilon_\mu(\nu).
    \end{alignat*}
    It follows from this and the $G$-invariance of $\langle -, -\rangle$ that
    \begin{alignat*}{4}
      \partial_{gs}(gf) &= g\partial_s(f) \qquad&&\text{and}\qquad&
        \partial_{gf}(gu) &= g\partial_f(u) \\
      \iota_{gm}(g\nu) &= g\iota_m(\nu)
        \qquad&&\text{and}\qquad&
       \iota_{g\mu}(g\ell) &= g\iota_\mu(\ell).
    \end{alignat*}
    The claimed $G$-equivariance follows.
  \end{proof}
\end{lemma}

\begin{example}
  A special case of the preceding construction gives
  a distinguished and familiar $G$-invariant endomorphism of
  $\rS(V^*) \otimes \wedge V^*$.
  Let $v_1, \ldots, v_n \in V$ be a basis and let
  $\lambda_1, \ldots, \lambda_n \in V^*$ be its dual basis. Then
    \[ \sum_{j=1}^n v_j \otimes \lambda_j \in \rS(V) \otimes \wedge V^* \]
  is independent of the choice of basis and is hence $G$-invariant.
  The action of this element is the
  \textit{exterior derivative}, namely
  \begin{equation}
    \dif \coloneqq \sum_{j=1}^n \partial_{v_j} \otimes \epsilon_{\lambda_j}
      \in \End_F(\rS(V^*) \otimes \wedge V^*).
  \end{equation}
  It satisfies $g(\dif \omega) = \dif(g\omega)$ for all $g \in G$,
  $\omega \in \rS(V^*) \otimes \wedge V^*$.
\end{example}

The operators above acting on $\rS(V^*) \otimes \wedge M^*$ generate the following
\textit{super Weyl algebra}. All of the operators that we will be using in this paper
come from the action of this algebra on $\rS(V^*) \otimes \wedge M^*$.

\begin{proposition}
  The subalgebra of $\End_F(\rS(V^*) \otimes \wedge M^*)$ generated by
  $m_\mu \otimes \id,
  \partial_v \otimes \id, \id \otimes \epsilon_\mu, \id \otimes \iota_v$ for
  $v \in V$, $\mu \in M^*$ is isomorphic to
  the tensor product of simple algebras
    \[ \bD(V) \otimes \Cliff(M \oplus M^*) \]
  where $\bD(V)$ is the algebra generated by $\partial_v, m_\lambda$ which
  is the Weyl algebra on $\dim(V)$ variables, and
  $\Cliff(M \oplus M^*)$ is the Clifford algebra of the split form
  $\langle m + \mu, \ell + \nu\rangle_{\mathrm{cl}} \coloneqq \mu(\ell) + \nu(m)$ which
  is generated by $\iota_m, \epsilon_\mu$.
\end{proposition}

\subsection{Derivations and anti-derivations}

The differential operators $\partial_s$ and $\partial_f$ are the usual
polynomial differential operators acting on the polynomial ring or its dual. The operators
$\epsilon_m$ and $\epsilon_\mu$ are called \textit{exterior products} and
the operators $\iota_m$ and $\iota_\mu$ are called \textit{interior products}. For
completeness and concreteness, we briefly summarize some of their well-known properties.

If $v \in V$, then $\partial_v \in \End_F(\rS(V^*))$ is given as usual by
\begin{equation}
  \partial_v f(x) = \left.\frac{\dif f(x+tv)}{\dif t}\right|_{t=0}.
\end{equation}
The operators $\partial_v$ satisfy the classical Leibniz rule
\begin{equation}
  \partial_v fh = (\partial_v f)h + f(\partial_v h)
\end{equation}
for all $v \in V$ and $f, h \in \rS(V^*)$, and hence are \textit{derivations}.
If $s \in \rS(V)$, we have
\begin{equation}
  (\partial_sf)(0) = \langle f, s\rangle = (\partial_fs)(0).
\end{equation}
Combining these last two observations, for all $v \in V$ and $\lambda \in V^*$, we have
\begin{equation}
  \partial_v m_\lambda - m_\lambda \partial_v = \lambda(v) \id_{\rS(V^*)}.
\end{equation}
It follows from the $\perm$ description in \Cref{lem:amp_pairings.1} that if
$v_1, \ldots, v_n \in V$ is a basis with dual basis $\lambda_1, \ldots, \lambda_n \in V^*$,
then
\begin{equation}
  \partial_{v^\alpha} \lambda^\beta
    = \delta_{\alpha \leq \beta} \frac{\beta!}{(\beta - \alpha)!} \lambda^{\beta - \alpha},
\end{equation}
where we have used multi-index notation and all operations are component-wise.

Analogously, one may check that if $m_1, \ldots, m_r \in M$ is a basis with dual basis
$\mu_1, \ldots, \mu_r \in M^*$, then
\begin{equation}
  \iota_{m^I} \mu^J = \pm \delta_{I \subset J} \mu^{J - I},
\end{equation}
where we have used the natural analogue of multi-index notation in this setting,
e.g.~$\mu^J \coloneqq \mu_{j_1} \wedge \cdots \wedge \mu_{j_\ell}$
if $J = \{j_1 < \cdots < j_\ell\}$. The sign may be
determined explicitly by iterating the well-known identity
\begin{equation}
  \iota_m(\mu_1 \wedge \cdots \wedge \mu_\ell)
    = \sum_{j=1}^\ell (-1)^{j-1} \mu_j(m)\, \mu_1 \wedge
        \cdots \wedge \widehat{\mu_j} \wedge \cdots \wedge \mu_\ell.
\end{equation}
More generally, the operators $\iota_m$ for $m \in M$ are
\textit{anti-derivations} in the sense that
\begin{equation}
  \iota_m(\mu \wedge \nu) = (\iota_m \mu) \wedge \nu + (-1)^k \mu \wedge (\iota_m \nu)
\end{equation}
for all $\mu \in \wedge^k M^*$ and $\nu \in \wedge M^*$. In particular, for all
$m \in M$ and $\xi \in M^*$, we have
\begin{equation}\label{eq:iota_eps}
  \iota_m \epsilon_\xi + \epsilon_\xi \iota_m = \xi(m) \id_{\wedge M^*}.
\end{equation}
See for instance \cite[pp.~356-359]{MR1930091} for further details.

\subsection{Hermitian forms}
While our key definitions and statements all involve the natural perfect pairings from
\Cref{ssec:pairings}, some of our proofs require replacing them with Hermitian forms.
For later use, we now recall some elementary properties of Hermitian forms and relate
the pairing $\rS(V^*) \otimes \wedge M^* \times \rS(V) \otimes \wedge M \to F$ to a
Hermitian form on $\rS(V^*) \otimes \wedge M^*$.
In this subsection, let $F$ be a subfield of $\bC$ closed under complex conjugation
and let $G \leq \GL(V)$ be an arbitrary \textit{finite} subgroup.

A \textit{Hermitian form} on an $F$-vector space $W$ is a map
  \[ (-, -) \colon W \times W \to F \]
such that if $w \in W$ is fixed, the map $v \mapsto (v, w)$
is $F$-linear, and $\overline{(v, w)} = (w, v)$. That is, a Hermitian form is linear
in the first argument and conjugate-linear in the second argument. Note that
$(w, w) = \overline{(w, w)} \in F \cap \bR$.
A Hermitian form is said to be \textit{positive-definite} if $(w, w) > 0$ for all
$0 \neq w \in W$.

\begin{lemma}\label{lem:V_Vperp}
  Let $W = \oplus_{j=0}^\infty W_j$ be a graded $F$-vector space with
  $\dim W_j < \infty$ for all $j$. Suppose that $(-, -)$ is a positive-definite
  Hermitian form on $W$ for which $(W_i, W_j) = 0$ if $i \neq j$. If $X$ is a
  graded subspace of $W$ and
    \[ X^\perp \coloneqq \{w \in W : (x, w) = 0, \forall x \in W\}, \]
  then $W = X \oplus X^\perp$.
  
  \begin{proof}
    It is enough to consider the case when $\dim W < \infty$. Note that if
    $x \in X \cap X^\perp$, then $(x, x) = 0$, so $x = 0$. We must show
    $X + X^\perp$ spans $W$. Let $x_1, \ldots, x_k$ be a basis for $X$.
    Then $X^\perp$ is the kernel of the $F$-linear map $W \to F^k$ given by
    $w \mapsto ((w, x_1), \ldots, (w, x_k))$. By the rank-nullity theorem,
    $\dim X + \dim X^\perp \geq \dim W$, and the result follows.
  \end{proof}
\end{lemma}

Since $F$ is closed under complex conjugation, every finite-dimensional $F$-vector
space $W$ has a positive-definite Hermitian form given by choosing a basis
$\omega_1, \ldots, \omega_k$ of $W^*$ and using
\begin{equation}
  (u, v) \coloneqq \sum_{i=1}^k \omega_i(u) \overline{\omega_i(v)}
\end{equation}
for all $u, v \in W$. We may strengthen this construction
and relate it to the canonical pairing
  \[ \langle -, -\rangle \colon \rS(V^*) \otimes \wedge M^* \times
      \rS(V) \otimes \wedge M \to F \]
from \Cref{ssec:pairings} as follows.

\begin{lemma}\label{lem:tau_herm}
  There is a positive-definite $G$-invariant Hermitian form
    \[ (-, -) \colon \rS(V^*) \otimes \wedge M^* \times \rS(V^*) \otimes \wedge M^* \to F \]
  and a conjugate-linear $G$-equivariant ring isomorphism
    \[ \tau \colon \rS(V) \otimes \wedge M \to \rS(V^*) \otimes \wedge M^* \]
  such that, for all $\eta \in \rS(V^*) \otimes \wedge M^*$ and
  $\omega \in \rS(V) \otimes \wedge M$,
  \begin{equation}\label{eq:tau_herm}
    \langle \eta, \omega\rangle = (\eta, \tau(\omega)).
  \end{equation}
  
  \begin{proof}
    Let $W$ be a finite-dimensional $G$-module.
    Since $F$ is assumed closed under complex conjugation, there is a positive-definite
    Hermitian form on $W$. Since $G$ is assumed finite, the Hermitian form may be taken
    to be $G$-invariant by Weyl's unitarian trick. We may extend $G$-invariant
    positive-definite Hermitian forms
    in three ways analogous to \Cref{lem:amp_pairings.1}.
    \begin{enumerate}[(1)]
      \item Suppose $W_1, W_2$ have $G$-invariant positive-definite
        Hermitian forms $(-, -)_i \colon W_i \times W_i \to F$. The form
        $W_1 \otimes W_2 \times W_1 \otimes W_2 \to F$ defined by
        $(w_1 \otimes w_2, w_1' \otimes w_2') \coloneqq (w_1, w_1')_1 (w_2, w_2')_2$
        extended bilinearly is a well-defined $G$-invariant Hermitian form. Moreover, it
        remains positive-definite as can be checked on orthogonal bases of $W_1, W_2$.
      \item A $G$-invariant positive-definite Hermitian form on $W$
        induces a $G$-invariant positive-definite Hermitian form on the
        $k$th exterior power $\wedge^k W$ by
        $(w_1 \wedge \cdots \wedge w_k, w_1' \wedge \cdots w_k')
        \coloneqq \det((w_i, w_j'))_{i, j=1}^k$.
      \item A $G$-invariant positive-definite Hermitian form on $W$
        induces a $G$-invariant positive-definite Hermitian form on the
        $k$th symmetric power $\rS^k(W)$ by
        $(w_1 \cdots w_k, w_1' \cdots w_k') \coloneqq \perm((w_i, w_j'))_{i, j=1}^k$,
        where $\perm$ denotes the matrix permanent.
    \end{enumerate}
    Combining these constructions yields a $G$-invariant positive-definite
    Hermitian form $(-, -)$ on $\rS(V^*) \otimes \wedge M^*$. We may
    thus define $\tau$ by \eqref{eq:tau_herm}.
    
    The conjugate-linearity and $G$-equivariance of $\tau$ follow quickly from
    \eqref{eq:tau_herm} and the corresponding properties of $(-, -)$ and
    $\langle -, -\rangle$. It remains to show that $\tau$ is multiplicative.
    On $\wedge M^*$, we have
    \begin{align*}
      (w_1 \wedge \cdots \wedge w_k, \tau(w_1' \wedge \cdots \wedge w_k'))
        &= \langle w_1 \wedge \cdots \wedge w_k, w_1' \wedge \cdots \wedge w_k'\rangle \\
        &= \det(\langle w_i, w_j'\rangle)
          = \det((w_i, \tau(w_j'))) \\
        &= (w_1 \wedge \cdots \wedge w_k, \tau(w_1') \wedge
          \cdots \wedge \tau(w_k')\rangle.
    \end{align*}
    The full calculation is exactly analogous.
  \end{proof}
\end{lemma}

%

\subsection{Pseudo-reflection groups and coinvariants}
We now recall some basic facts concerning pseudo-reflection groups and establish
some further notation. We again let $F$ be an arbitrary field of characteristic $0$.

\begin{definition}
  A \textit{pseudo-reflection} is an element $g \in \GL(V)$ such that
    \[ \dim \ker(g-I) = n-1. \]
  A pseudo-reflection of order two is called a \textit{reflection}. A
  \textit{(pseudo) reflection group} is a finite subgroup of $\GL(V)$ generated
  by (pseudo) reflections.
\end{definition}

\noindent For the rest of this subsection, let $G$ denote a pseudo-reflection group.

Shephard--Todd \cite{MR0059914} and Chevalley \cite{MR72877}
showed that the $G$-invariants $\rS(V^*)^G$ are generated by $n$
algebraically independent, homogeneous elements $f_1, \ldots, f_n \in \rS(V^*)^G$
of positive degrees $d_1, \ldots, d_n$.  The $f_i$ are called \textit{basic invariants}
of $G$ and the $d_i$ are called the \textit{degrees} of $G$.

Recall that the \textit{Hilbert series} of a graded
vector space $M = \oplus_{i=0}^\infty M_i$ with $\dim M_i < \infty$
is the formal power series
  \[ \Hilb(M; q) \coloneqq \sum_{i=0}^\infty q^i \dim M_i. \]
When $M = \oplus_{i,j=0}^\infty M_{i,j}$ is bigraded, we use a bivariate power series
  \[ \Hilb(M; q, z) \coloneqq \sum_{i,j=0}^\infty q^i z^j \dim M_{i,j}. \]
  
If $g \in \GL(V)$ is a
(pseudo) reflection, then $g^* \in \GL(V^*)$ is as well, so $G^* \leq \GL(V^*)$
is a pseudo-reflection group, and so $\rS(V)^G$ is similarly
generated by algebraically independent, homogeneous elements
$z_1, \ldots, z_n \in \rS(V)^G$. Since $G$ acts completely reducibly,
$\dim M^G = \dim (M^*)^G$ in general, so
$\Hilb(\rS(V)^G; q) = \Hilb(\rS(V^*)^G; q)$, and we may take
$\deg f_i = \deg z_i$.

\begin{definition}
  The \textit{coinvariant ideal} of $\rS(V)$ is the ideal
  $\cI$ generated by all homogeneous $G$-invariants of positive degree.
  We write $\cI^*$ for the coinvariant ideal of $\rS(V^*)$. The
  \textit{coinvariant algebra} of $\rS(V)$ or $\rS(V^*)$ is
  $\rS(V)/\cI$ or $\rS(V^*)/\cI^*$, respectively.
\end{definition}

\begin{remark}\label{rem:chev}
  Chevalley \cite{MR72877} showed that $\rS(V^*)/\cI^*$, a graded
  algebra, carries the regular representation of $G$ and that
  \begin{equation}\label{chev_tensor}
    \rS(V^*)^G \otimes \rS(V^*)/\cI^* \cong \rS(V^*)
  \end{equation}
  as a graded $\rS(V^*)^G$-module.
\end{remark}

\subsection{Fields of definition}
In order to use the machinery of Hermitian forms, we require representations
defined over subfields of the complex numbers which are closed under complex conjugation.
In practice, pseudo-reflection groups are typically constructed in terms of explicit
unitary matrices over cyclotomic fields, in which case these properties are trivial.
More care is required to handle the general case and avoid artificial assumptions.

Suppose $W$ is an $n$-dimensional vector
space over a field $K$ of characteristic $0$. The \textit{character field}
of a representation $\rho \colon G \to \GL(W)$ is the subfield
$\bQ(\rho) \subset K$ generated by the set $\{\Tr\rho(g) \mid g \in G\}$.

\begin{lemma}
  If $|G| < \infty$, then $\bQ(\rho)$ is isomorphic with a subfield of $\bC$
  closed under complex conjugation.
  
  \begin{proof}
    Choosing some basis for $W$, the matrix coefficients of $\rho(g)$ are
    $c_{ij}^g \in K$. Clearly $\bQ(\rho) \subset \bQ(c_{ii}^g)$. Since
    $\{c_{ii}^g\}$ is finite, we may identify $\bQ(c_{ii}^g)$ with a subfield of $\bC$, so
    $\bQ(\rho)$ may be identified with a subfield of $\bC$.
    Diagonalizing now shows that $\Tr\rho(g)$ is of the form
    $\sum_i \zeta^{\alpha_i}$ for $\zeta = \exp(2\pi i/|G|) \in \bC$. Thus 
      \[ \Tr\rho(g^{-1}) = \sum_i \zeta^{-\alpha_i} = \sum_i \overline{\zeta^{\alpha_i}}
          = \overline{\Tr\rho(g)}, \]
    so $\bQ(\rho) \subset \bC$ is closed under complex conjugation.
  \end{proof}
\end{lemma}

We say $\rho \colon G \to \GL(W)$ is \textit{defined over a subfield} $K'$ of $K$ if
there is some basis of $W$ for which $\rho(G) \subset \GL_n(K')$.
If $\rho$ is defined over $K'$, then clearly $K' \supset \bQ(\rho)$.
In favorable circumstances, the representation is actually defined over its character field.

\begin{theorem}[{Clark--Ewing \cite{MR367979}; see \cite[Appendix~B, p.~359]{MR1838580}}]\label{thm:Clark-Ewing}
  If $\rho \colon G \to \GL(W)$ is an irreducible representation of a finite
  group $G$ over a field $K$ of characteristic $0$ and $\rho(G)$
  contains a pseudo-reflection, then $\rho$ is defined over $\bQ(\rho)$.
\end{theorem}

Benard gave a similar result for all representations of
pseudo-reflection groups in characteristic $0$. The only known
proofs are case-by-case using the Shephard--Todd classification.

\begin{theorem}[{Benard \cite[Thm.~1]{MR401901}}]\label{thm:benard}
  Let $G$ be a pseudo-reflection group over a field of characteristic $0$. Then
  the character field $K$ of $G$ is a splitting field of $G$. That is, every representation
  of $G$ over a field containing $K$ is defined over $K$.
\end{theorem}

\begin{corollary}\label{cor:complex_conj}
  If $G$ is a pseudo-reflection group over a field $K$ of characteristic $0$ and
  $\rho \colon G \to \GL(W)$ is a representation of $G$ over $K$, then $\rho$
  is defined over a subfield of $K$ which is isomorphic to a subfield of $\bC$
  closed under complex conjugation.
\end{corollary}

\section{Exponents, Jacobians, and Vandermondians}\label{sec:exps}

We continue to let $V$ be an $n$-dimensional vector space over a field $F$ of
characteristic $0$. Let $G \leq \GL(V)$ be a pseudo-reflection group and
suppose $M$ is an $r$-dimensional $G$-module. In this section we recall the
$M$-exponents, the basic derivations for $M$, and related notions.

\subsection{Exponents}
Chevalley's results in \Cref{rem:chev} (cf.~\cite[Lemma~6.45]{MR1217488})
imply that there are homogeneous elements
$\omega_1^M, \ldots, \omega_r^M \in (\rS(V^*) \otimes M^*)^G$ of
bi-degrees $(e_1^M, 1), \ldots, (e_r^M, 1)$ such that
\begin{equation}\label{eq:omegas}
  (\rS(V^*) \otimes M^*)^G = \rS(V^*)^G\omega_1^M \oplus \cdots
      \oplus \rS(V^*)^G\omega_r^M.
\end{equation}

\begin{definition}\label{def:exponents}
  The $\omega_i^M$ are called \textit{basic derivations} for $M$ over $\rS(V^*)$
  (cf.~\cite[Def.~6.50]{MR1217488}). The $e_i^M$ are the \textit{$M$-exponents}.
  The \textit{exponents of $G$} are the $V$-exponents
  $e_1, \ldots, e_n \coloneqq e_1^V, \ldots, e_n^V$,
  and the \textit{coexponents of $G$} are the $V^*$-exponents
  $e_1^*, \ldots, e_n^* \coloneqq e_1^{V^*}, \ldots, e_n^{V^*}$.
\end{definition}

The basic derivations are not unique, but the $M$-exponents are uniquely determined
by $M$ up to rearrangement. If $M$ is absolutely irreducible, the $M$-exponents
are the degrees in which $M$ appears in $\rS(V^*)/\cI^*$.

Similarly, the \textit{basic derivations} of a $G$-module $M$ of
dimension $r$ over $\rS(V)$ (rather than $\rS(V^*)$) are homogeneous elements
$\tilde{\omega}_1^M, \ldots, \tilde{\omega}_r^M \in (\rS(V) \otimes M)^G$ such that
\begin{equation}
  (\rS(V) \otimes M)^G = \rS(V)^G\tilde{\omega}_1^M \oplus \cdots
    \oplus \rS(V)^G\tilde{\omega}_r^M.
\end{equation}
After rearrangement, $\bideg \tilde{\omega}_i^M = \bideg \omega_i^M$ since
$\Hilb((\rS(V) \otimes M)^G; q) = \Hilb((\rS(V^*) \otimes M^*)^G; q)$,

\begin{example}\label{ex:derivations_exterior}
  Since basic invariants $f_1, \ldots, f_n \in \rS(V^*)^G$ are
  algebraically independent and $\operatorname{char} F = 0$, the classical Jacobian criterion
  gives $\det(\partial_{v_j} f_i) \neq 0$ for any basis $v_1, \ldots, v_n$ of $V$.
  This observation and the well-known fact that $e_i = d_i - 1$ imply
  $\dif f_1, \ldots, \dif f_n \in (\rS(V^*) \otimes V^*)^G$
  form a set of basic derivations for $V$.
  That is, we may take
    \[ \omega_i^V = \dif f_i = \sum_{j=1}^n \partial_{v_j} f_i \otimes \lambda_j, \]
  where $\lambda_1, \ldots, \lambda_n$ is the dual basis of $v_1, \ldots, v_n$.

  Likewise, the basic invariants $z_1, \ldots, z_n \in \rS(V)^G$ yield basic
  derivations $\tilde{\omega}_i^V = \dif z_i
  = \sum_{j=1}^n \partial_{\lambda_j} z_i \otimes v_j \in (\rS(V) \otimes V)^G$.
\end{example}

\subsection{Jacobians and Vandermondians}
The preceding example motivates the following general notion.
Suppose $M^*$ has basis $\mu_1, \ldots, \mu_r$. We may expand the basic
derivations for $M$ over $\rS(V^*)$ as
  \[ \sum_{j=1}^r J_{ij}^M \otimes \mu_j \coloneqq \omega_i ^M
      \in (\rS(V^*) \otimes \wedge M^*)^G \]
for some $J_{ij}^M \in \rS(V^*)$. 

\begin{definition}\label{def:jacobian}
  The \textit{Jacobian} of $M$ is
    \[ J_M \coloneqq \det(J_{ij}^M)_{i,j=1}^r \in \rS(V^*). \]
\end{definition}

The Jacobian is non-zero and is uniquely determined up to a
non-zero constant \cite[Prop.~2.5]{MR575083}. Note that
$\deg J_M = e_1^M + \cdots + e_r^M$.

\begin{example}
  We have
    \[ \omega_1^V \wedge \cdots \wedge \omega_n^V
        = \dif f_1 \wedge \cdots \wedge \dif f_n
        = J_V\,\eta
            \in \rS(V^*) \otimes \wedge^n V^* \]
  for some $0 \neq \eta \in \wedge^n V^*$.
  Since this expression is $G$-invariant, it follows that $J_V$ belongs to
  $\rS(V^*)^{\det_V}$, the $\det_V$-isotypic component of $\rS(V^*)$.
  Steinberg \cite{MR117285} proved that
  \begin{equation}\label{eq:steinberg_JV}
    \rS(V^*)^{\det_V} = \rS(V^*)^G J_V.
  \end{equation}
  Indeed, $e_i = d_i - 1$ follows from \eqref{eq:steinberg_JV}.
  For general $M$, we similarly have some $0 \neq \xi \in \wedge^r M$ such that
  \begin{equation}
      \omega_1^M \wedge \cdots \wedge \omega_r^M
      = J_M\,\xi
      \qquad\text{and}\qquad J_M \in \rS(V^*)^{\det_M}.
  \end{equation}
  Consequently,
  \begin{equation}
    \rS(V^*)^{\det_M} \supset \rS(V^*)^G J_M,
  \end{equation}
  though the containment may be strict, which
  motivates the following.
\end{example}

By \eqref{eq:omegas}, there
exists an element
  \[ \Delta_M \in \rS(V^*)^{\det_M} \]
such that
  \[ \rS(V^*)^{\det_M} = \rS(V^*)^G \Delta_M. \]
\begin{definition}\label{def:Vandermondian}
  We call $\Delta_M$ the \textit{Vandermondian of $M$ over $\rS(V^*)$}.
\end{definition}

The Vandermondian of $M$ is non-zero and is uniquely determined up to a non-zero scalar.
When $G$ consists of permutation matrices, $\Delta_V$ is the classical Vandermondian
determinant. Steinberg's result \eqref{eq:steinberg_JV} implies
  \[ \Delta_V = J_V \]
up to a non-zero scalar. In general, only
  \[ \Delta_M \mid J_M \]
is true. We have $\Delta_M = J_{\det M} = \Delta_{\det M}$ up to a non-zero scalar,
as can be seen from \eqref{eq:omegas}. In particular, when
$M$ is one-dimensional, we may take $\Delta_M = J_M$.

Interchanging $V$ and $V^*$, we have the Jacobian and Vandermondian
with respect to $\rS(V)$ instead of $\rS(V^*)$.

\begin{definition}
  The \textit{Jacobian of $M$ with respect to $\rS(V)$} is defined by
    \[ \tilde{J}_M\,\eta = \tilde{\omega}_1^M \wedge \cdots \wedge \tilde{\omega}_r^M
        \in (\rS(V) \otimes \wedge M)^G \]
  for $0 \neq \eta \in \wedge M^r$.
  The \textit{Vandermondian of $M$ with respect to $\rS(V)$} is defined by
    \[ \rS(V)^{\det_{M^*}} = \rS(V)^G \tilde{\Delta}_M. \]
\end{definition}

As before, $\tilde{J}_M$ and $\tilde{\Delta}_M$ are each determined by $M$ up to
a non-zero scalar. Since $\bideg \tilde{\omega}_i^M = \bideg \omega_i^M$, we have
$\deg \tilde{J}_M = \deg J_M$ and $\deg \tilde{\Delta}_M = \deg \Delta_M$.

\subsection{Gutkin's formula}
Stanley \cite{MR460484} expressed $\Delta_M$ and Gutkin
\cite{MR0314956} expressed $J_M$ as a product of
linear forms vanishing on the reflecting hyperplanes of $G$, which we now summarize.
See \cite[\S4]{MR3959046} for historical discussion and
\cite[\S4.5.2]{MR2590895} for a proof using Molien's theorem.

Let $\cA(G)$ be the set of reflecting hyperplanes of $G$, i.e.~the fixed spaces of
pseudo-reflections of $G$. For each $H \in \cA(G)$, fix some $\alpha_H \in V^*$
with $\ker \alpha_H = H$. Let $G_H$ denote the subgroup of
$G$ fixing $H$ pointwise. Since $\sigma \in G_H$ acts trivially on a codimension $1$
subspace, $\sigma \mapsto \det_{V^*} \sigma$ is a faithful representation
$G_H \hookrightarrow F^\times$. Since finite subgroups of $F^\times$ are cyclic,
$G_H$ is cyclic, $G_H$ is generated by a pseudo-reflection, and the irreducible
$G_H$-representations are powers of $\det_{V^*}$. Consequently,
\begin{equation}\label{eq:M_res}
  M\res_{G_H}^G \cong {\det}_{V^*}^{m_{H, 1}(M)} \oplus \cdots
    \oplus {\det}_{V^*}^{m_{H, r}(M)}
\end{equation}
as $G_H$-modules, with $0 \leq m_{H, i}(M) < |G_H|$. Set
  \[ m_H(M) \coloneqq m_{H, 1}(M) + \cdots + m_{H, r}(M). \]

\begin{theorem}[{Gutkin \cite{MR0314956}; cf. \cite[(2.11)]{MR575083} or
\cite[Theorem~4.38(1)]{MR2590895}}]\label{thm:Gutkin}
  Up to a non-zero scalar,
    \[ J_M = \prod_{H \in \cA(G)} \alpha_H^{m_H(M)}. \]
\end{theorem}

Recall that $\Delta_M = J_{\det_M}$ up to a non-zero scalar, so \Cref{thm:Gutkin}
gives product formulas for $\Delta_M$ as well.
We have the following well-known special cases. See \cite[Prop.~4.34(2)]{MR2590895}
and \cite{MR575083}.

\begin{corollary}\label{cor:Gutkin}
  \ 
  \begin{enumerate}[(a)]
    \item $J_{V^*} = \prod_{H \in \cA(G)} \alpha_H$.
    \item $J_V = \prod_{H \in \cA(G)} \alpha_H^{|G_H| - 1}$.
    \item $J_M$ equals $\Delta_M$ up to a non-zero scalar if and only if $m_H(M) < |G_H|$ for
      all $H \in \cA(G)$.
    \item\label{cor:Gutkin.pseudos}
      If the pseudo-reflections of $G$ act on $M$ as pseudo-reflections or the identity,
      then $J_M$ equals $\Delta_M$ up to a non-zero scalar.
    \item If $M$ is a Galois conjugate of $V$, then $J_M$ equals $\Delta_M$ up to a
      non-zero scalar.
    \item If $\chi$ is one-dimensional, then $\Delta_\chi \mid \Delta_V$.
  \end{enumerate}
  
  \begin{proof}
    For (a), use $V^*\res_{G_H}^G \cong \det_{V^*} \oplus 1^{n-1}$. For (b), use
    $V\res_{G_H}^G \cong \det_V = {\det}_{V^*}^{-1} \oplus 1^{n-1}$. For (c),
    \eqref{eq:M_res} implies
      \[ m_H(M) - m_H({\det}_M) \in \bZ_{\geq 0} |G_H|, \]
    which gives the result since $\Delta_M = J_{\det_M}$.
    Now (d) follows from (c) since in this case
    \eqref{eq:M_res} has at most one non-trivial summand, so $m_H(M) < |G_H|$,
    and (e) is a special case of (d). Finally, (f) follows from (b) and (d) since
    $m_H(\chi) = m_{H, 1}(\chi) < |G_H|$ and $\Delta_V = J_V$.
  \end{proof}
\end{corollary}

\section{Harmonics and semi-invariant bases}\label{sec:super_alternants}
Our first goal in this section is to define the well-known $G$-harmonic polynomials using
the constructions in \Cref{sec:pairings}. We
prove their basic properties by constructing a $G$-invariant Hermitian
form, and then use the harmonics to construct explicit bases for semi-invariant differential
forms. The main results of this section are
\Cref{thm:super_alternant_basis} and \Cref{thm:super_alternant_basis.2}.
$F$ remains a field of characteristic $0$, $V$ is an $n$-dimensional
$F$-vector space, $G \leq \GL(V)$ is a
pseudo-reflection group, and $M$ is an $r$-dimensional $G$-module.

\begin{definition}\label{def:harmonics}
  The spaces of \textit{$G$-harmonic} elements of $\rS(V^*)$ and $\rS(V)$ are,
  respectively,
   \begin{align*}
     \cH(V^*)
       &\coloneqq \{f \in \rS(V^*) : \langle s, f\rangle = 0 \text{ for all }s \in \cI\}, \\
     \cH(V)
       &\coloneqq \{s \in \rS(V) : \langle s, f\rangle = 0 \text{ for all }f \in \cI^*\}.
   \end{align*}
\end{definition}

That is, the harmonics are the orthogonal complements
of the coinvariant ideals $\cI$ and $\cI^*$ with respect to the natural perfect pairing
$\langle -, -\rangle$ from \eqref{eq:perfect_pairing}. The harmonics have the
following basic properties.

\begin{lemma}\label{lem:harmonics_props.1}
  We have
  \begin{align}
    \label{eq:harmonics.1}
    \cH(V^*) &= \{f \in \rS(V^*) : \partial_{z_j} f = 0, j=1, \ldots, n\}, \\
    \label{eq:harmonics.2}
    \rS(V^*) &= \cH(V^*) \oplus \cI^*, \\
    \label{eq:harmonics.3}
    \cH(V^*) &\cong \rS(V^*)/\cI^* \quad\text{as graded $G$-modules}, \\
    \label{eq:harmonics.4}
    \rS(V^*) &= \rS(V^*)^G \otimes \cH(V^*),
  \end{align}
  and similarly with $\cH(V)$.
  
  \begin{proof}
    \eqref{eq:harmonics.1} follows from the fact that $\cI$ is the ideal
    generated by $z_1, \ldots, z_n$. \eqref{eq:harmonics.4} follows from
    \eqref{eq:harmonics.3} and Chevalley's result in \Cref{rem:chev}, and
    \eqref{eq:harmonics.3} follows from \eqref{eq:harmonics.2}.
    
    As for \eqref{eq:harmonics.2}, first suppose $F \subset \bC$ is
    closed under complex conjugation. By \Cref{lem:tau_herm}, we have
    $\tau(\rS(V)^G) = \rS(V^*)^G$ and $\tau(\cI) = \cI^*$.
    Thus $\langle \cI, f\rangle = (\tau(\cI), f) = (\cI^*, f)$,
    so $\cH(V^*)$ is the orthogonal complement of $\cI^*$ under
    a positive-definite Hermitian form. Now \eqref{eq:harmonics.2} follows
    from \Cref{lem:V_Vperp} in this case.
    
    For the general case, by \Cref{cor:complex_conj}, $G$ is defined over a subfield
    $K$ of $F$ which may be identified with a subfield of $\bC$ closed under
    complex conjugation. Let $v_1, \ldots, v_n$ be a basis for $V$ over $F$ such that the
    matrix of each $\sigma \in G$ is in $\GL_n(K)$. Set $V_K \coloneqq \oplus_{i=1}^n Kv_i$.
    Then
      \[ \rS(V^*) = \Span_F \rS(V_K^*) \cong F \otimes_K \rS(V_K^*) \]
    as $F$-algebras and graded $G$-modules.
    Let $R_G \coloneqq \frac{1}{|G|} \sum_{\sigma \in G} \sigma$
    be the projection onto $G$-invariants. Since $R_G$ is $F$-linear,
    \begin{align*}
      \rS(V^*)^G
        &= R_G \rS(V^*)
        \cong R_G(F \otimes_K \rS(V_K^*)) \\
        &= F \otimes_K R_G \rS(V_K^*)
        = F \otimes_K \rS(V_K^*)^G.
    \end{align*}
    This implies $\cI^* \cong F \otimes_K \cI_K^*$, and hence
    $\cH(V^*) \cong F \otimes_K \cH(V_K^*)$.
  \end{proof}
\end{lemma}

\begin{remark}
  By \Cref{lem:tau_herm}, when $F \subset \bC$ is closed under complex conjugation,
  we have $\tau((\rS(V) \otimes M)^G)
  = (\rS(V^*) \otimes M^*)^G$, so we may take
  $\tau(\tilde{\omega}_i^M) = \omega_i^M$, $\tau(\tilde{J}_M) = J_M$, and
  $\tau(\tilde{\Delta}_M) = \Delta_M$.
\end{remark}

Recall the $G$-equivariant $F$-algebra homomorphism 
\begin{align*}
  \rS(V) \otimes \wedge M^* &\to \End_F(\rS(V^*) \otimes \wedge M^*) \\
  s \otimes \mu &\mapsto \partial_s \otimes \epsilon_\mu.
\end{align*}
from \Cref{lem:actions}.

\begin{definition}\label{def:d_i_star}
  For $i=1, \ldots, r$, let
    \[ d_i^* \colon \rS(V^*) \otimes \wedge M^* \to \rS(V^*) \otimes \wedge M^* \]
  denote the action of the basic derivation
  $\tilde{\omega}_i^{M^*} \in (\rS(V) \otimes \wedge M^*)^G$
  from \Cref{sec:exps}.
\end{definition}

\begin{theorem}\label{thm:super_alternant_basis}
  Suppose $J_{M^*} \mid \Delta_\chi$. Then the $2^r$ elements
    \[ \{d_{i_1}^* \cdots d_{i_k}^* \Delta_{\chi}
      \mid 1 \leq i_1 < \cdots < i_k \leq r\} \]
  form a basis of $(\cH(V^*) \otimes \wedge M^*)^{\chi}$, and
  hence descend to a basis of $(\rS(V^*)/\cI^* \otimes \wedge M^*)^{\chi}$.
  
  \begin{proof}
    We suppose throughout that $F$ is a subfield of $\bC$ closed under
    complex conjugation. The general case follows by arguing
    as in the proof of \Cref{lem:harmonics_props.1} using extension of scalars.
    
    We have $\Delta_\chi \in \cH(V^*)^\chi$ since $\Delta_\chi$ is the lowest-degree
    element of $\rS(V^*)^\chi$ and so $\partial_{z_j} \Delta_\chi = 0$ for $j=1, \ldots, n$.
    Since $\tilde{\omega}_i^{M^*}$ is $G$-invariant,
    $d_i^*$ preserves the $\chi$-isotypic component.
    Since $\partial_s$ preserves $\cH(V^*)$, $\partial_s \otimes \epsilon_\mu$ preserves
    $\cH(V^*) \otimes \wedge M^*$, so $d_i^*$ preserves $\cH(V^*) \otimes \wedge M^*$
    and $d_{i_1}^* \cdots d_{i_k}^* \Delta_\chi \in (\cH(V^*) \otimes \wedge M^*)^\chi$.
    
    We now prove linear independence. Suppose
    \begin{equation}\label{eq:zeta.1}
      \sum \zeta_{i_1, \ldots, i_k} d_{i_1}^* \cdots d_{i_k}^* \Delta_{\chi} = 0
    \end{equation}
    for scalars $\zeta_{i_1, \ldots, i_k} \in F$. By homogeneity, we may suppose
    the subsets are each of size $k$. For a fixed choice of $\{i_1, \ldots, i_k\}$,
    let $\{j_1, \ldots, j_{r-k}\}$ be its complement in $[r]$. Since
    $\tilde{\omega}_i^{M^*} \tilde{\omega}_j^{M^*} = -\tilde{\omega}_j^{M^*}
    \tilde{\omega}_i^{M^*}$,
    we have $d_i^* d_j^* = -d_j^* d_i^*$. Applying $d_{j_1}^* \cdots d_{j_{r-k}}^*$
    to \eqref{eq:zeta.1} gives
    \begin{equation}\label{eq:zeta.2}
      \pm \zeta_{i_1, \ldots, i_k} d_1^* \cdots d_r^* \Delta_{\chi} = 0.
    \end{equation}
    We have $\tilde{\omega}_1^{M^*} \cdots \tilde{\omega}_r^{M^*}
    = \tilde{J}_{M^*} \eta$ where $0 \neq \eta \in \wedge^r M^*$, so
    $d_1^* \cdots d_r^* = \partial_{\tilde{J}_{M^*}} \otimes \epsilon_\eta$. Thus
    \eqref{eq:zeta.2} becomes
    \begin{equation}\label{eq:zeta.3}
      \zeta_{i_1, \ldots, i_k} \partial_{\tilde{J}_{M^*}} \Delta_{\chi} \eta = 0.
    \end{equation}
    Since $\Delta_\chi/J_{M^*} \in \rS(V^*)$ by assumption,
    $\tilde{\Delta}_\chi/\tilde{J}_{M^*} = \tau^{-1}(\Delta_\chi/J_{M^*}) \in \rS(V)$.
    Applying $\partial_{\tilde{\Delta}_\chi/\tilde{J}_{M^*}}$ to \eqref{eq:zeta.3}
    and canceling $\eta$ gives
    \begin{align*}
      0 &= \zeta_{i_1, \ldots, i_k} \partial_{\tilde{\Delta}_\chi} \Delta_\chi
          = \zeta_{i_1, \ldots, i_k} \langle \tilde{\Delta}_\chi, \Delta_\chi\rangle \\
        &= \zeta_{i_1, \ldots, i_k} (\tau(\tilde{\Delta}_\chi), \Delta_\chi)
          = \zeta_{i_1, \ldots, i_k} (\Delta_\chi, \Delta_\chi).
    \end{align*}
    Since $(\Delta_\chi, \Delta_\chi) > 0$, we find $\zeta_{i_1, \ldots, i_k} = 0$.
    
    From Chevalley's result in \Cref{rem:chev}, $\cH(V^*)$ carries the regular
    representation $F[G]$ of $G$. It is well-known that $F[G] \otimes M$ is
    isomorphic to $\dim M$ copies of $F[G]$ as a $G$-module, since if $M'$ is $M$
    with trivial $G$-action,
    \begin{align*}
      F[G] \otimes M &\to F[G] \otimes M' \\
      g \otimes m &\mapsto g \otimes g^{-1}m
    \end{align*}
    is a $G$-module isomorphism. Thus
    $\dim (\cH(V^*) \otimes \wedge M^*)^\chi = \dim \wedge M^* = 2^r$,
    so the linearly independent set is indeed a basis.
  \end{proof}
\end{theorem}

\begin{remark}
  The $d_i^*$ essentially appear in \cite[\S6, Prop.~18]{MR2163712} in slightly less
  generality. The preceding proof exploits their exterior algebra structure in a fundamental way.
  The condition $J_{M^*} \mid \Delta_\chi$ is a generalization of the condition from
  \cite{MR2163712} that $\chi$ is \textit{wholly-nontrivial}, meaning that
  $\Delta_{\det_{V^*}} \mid \Delta_\chi$.
\end{remark}

\begin{corollary}\label{cor:super_alternant_basis}
  Suppose $J_{M^*} \mid \Delta_\chi$. Then
    \[ \Hilb((\cH(V^*) \otimes \wedge M^*)^\chi; q, z)
        = q^{\deg \Delta_\chi} \prod_{i=1}^r (1 + zq^{-e_i^{M^*}}). \]
  
  \begin{proof}
    The $d_i^*$ alter bidegree by $(-e_i^{M^*}, 1)$.
  \end{proof}
\end{corollary}

We additionally have the following result. Alternatively, it is a straightforward
consequence of Orlik--Solomon's generalization \cite[Thm.~3.1]{MR575083}
of Solomon's description \cite{MR0154929} of $(\rS(V^*) \otimes \wedge V^*)^G$
as the exterior algebra over $\rS(V^*)^G$ generated by $\dif f_i$. By
\Cref{cor:Gutkin}(\ref{cor:Gutkin.pseudos}), $J_M = \Delta_M$ up to a non-zero scalar
whenever the pseudo-reflections of $G$ act on $M$ as pseudo-reflections or the identity.

\begin{corollary}
  Suppose $J_M = \Delta_M$ up to a non-zero scalar. Then
  \begin{equation}\label{eq:wedge_mult}
    \Hilb(\Hom_G(\wedge^k M, \rS(V^*)/\cI^*); q)
      = \sigma_k(q^{e_1^M}, \ldots, q^{e_r^M}),
  \end{equation}
  where $\sigma_k$ denotes an elementary symmetric polynomial of
  degree $k$. In particular, if $\wedge^k M$ is absolutely irreducible, 
  the right-hand side of \eqref{eq:wedge_mult} is the generating function
  for the degrees in which $\wedge^k M$ appears in the coinvariant
  algebra $\rS(V^*)/\cI^*$.
  
  \begin{proof}
    We have
      \[ (\rS(V^*)/\cI^* \otimes \wedge^k M^*)^G \cong
          \Hom_G(\wedge^k M, \rS(V^*)/\cI^*), \]
    so \eqref{eq:wedge_mult} is equivalent to
      \[ \Hilb((\rS(V^*)/\cI^* \otimes \wedge M^*)^G; q, z)
          = \prod_{i=1}^r (1 + zq^{e_i^M}). \]
    Since $\deg \Delta_M = \deg J_M = \sum_{i=1}^r e_i^M$,
    \Cref{cor:super_alternant_basis} gives
      \[ \Hilb((\rS(V^*)/\cI^* \otimes \wedge M)^{\det_M}; q, z)
          = \prod_{i=1}^r (z + q^{e_i^M}) = z^r \prod_{i=1}^r (1 + z^{-1} q^{e_i^M}). \]
    It thus suffices to show
      \[ \Hilb((\rS(V^*)/\cI^* \otimes \wedge M)^{\det_M}; q, z)
          = z^r \Hilb((\rS(V^*)/\cI^* \otimes \wedge M^*)^G; q, z^{-1}). \]
    This follows from the well-known fact that for each $0 \leq k \leq r$,
    there is a (non-canonical) isomorphism of $\GL(V)$-modules
    \begin{equation*}
      \wedge^k M^* \cong {\det}_M^{-1} \otimes \wedge^{r-k} M.
    \end{equation*}
  \end{proof}
\end{corollary}

\begin{remark}
  By Benard's \Cref{thm:benard}, every irreducible representation of
  a pseudo-reflection group $G$ is absolutely irreducible over the character field of $G$.
  Steinberg noted that if $G$ is generated
  by $n=\dim(V)$ pseudo-reflections and $V$ is irreducible, then $\wedge^k V$ is
  (absolutely) irreducible for all $0 \leq k \leq n$
  (see \cite[Ch.~V, \S2, Exercise 3(d)]{MR0240238}, \cite[p.~250, Thm.~A]{MR1838580}).
\end{remark}

We now restate and prove \Cref{thm:super_alternant_basis.2} and
\Cref{cor:super_alternant_basis.2} from the Introduction.

\begin{theorem}\label{thm:super_alternant_basis.2}
  Suppose $J_{M^*} \mid \Delta_\chi$. Then either of the sets
    \[ \{f_1^{a_1} \cdots f_n^{a_n} d_{i_1}^* \cdots d_{i_k}^* \Delta_\chi
           : 1 \leq i_1 < \cdots < i_k \leq r, a_j \in \bZ_{\geq 0}\} \]
  or
    \[ \{d_{i_1}^* \cdots d_{i_k}^* f_1^{a_1} \cdots f_n^{a_n} \Delta_\chi
           : 1 \leq i_1 < \cdots < i_k \leq r, a_j \in \bZ_{\geq 0}\} \]
  form bases for $(\rS(V^*) \otimes \wedge M^*)^\chi$.
    
  \begin{proof}
    The fact that the first set is a basis follows from
    \Cref{thm:super_alternant_basis} and \eqref{eq:harmonics.4}. For the second set,
    it suffices to verify linear independence. To do so, we sketch how to modify the
    proof of \Cref{thm:super_alternant_basis}. Suppose $F$ is a subfield of $\bC$
    closed under complex conjugation, pick a homogeneous, orthogonal basis
    $\{\overline{g}_\alpha\}$ of $\rS(V^*)^\chi = \rS(V^*)^G \Delta_\chi$
    with respect to the Hermitian form, and let $g_\alpha \Delta_\chi = \overline{g}_\alpha$
    where $\{g_\alpha\}$ is a basis for $\rS(V^*)^G$. Linear independence is
    unaffected if we replace $\{f_1^{a_1} \cdots f_n^{a_n}\}$ with $\{g_\alpha\}$.
    If
      \[ 0 = \sum_{\alpha, I} \zeta_{\alpha; I} d_{i_1}^* \cdots d_{i_k}^*
          g_\alpha \Delta_\chi, \]
    we may apply $\partial_{\tilde{\Delta}_\chi/\tilde{J}_{M^*}}
    \partial_{\tau^{-1}(g_\beta)} d_{j_1}^* \cdots d_{j_{r-k}}^*$
    to get
    \begin{align*}
      0 &= \sum_\alpha \pm \zeta_{\alpha; I} \partial_{\tilde{\Delta}_\chi/\tilde{J}_{M^*}}
        \partial_{\tau^{-1}(g_\beta)} \partial_{\tilde{J}_{M^*}}
        g_\alpha \Delta_\chi \eta \\
         &= \sum_\alpha \pm \zeta_{\alpha; I}
           \langle\tau^{-1}(g_\beta \Delta_\chi), 
           g_\alpha \Delta_\chi\rangle \eta \\
        &= \sum_\alpha \pm \zeta_{\alpha; I} (\overline{g}_\beta, \overline{g}_\alpha) \eta
          = (\overline{g}_\beta, \overline{g}_\beta) \zeta_{\beta; I} \eta.
     \end{align*}
     Now use $(\overline{g}_\beta, \overline{g}_\beta) > 0$.
  \end{proof}
\end{theorem}

\begin{corollary}\label{cor:super_alternant_basis.2}
  Suppose $J_{M^*} \mid \Delta_\chi$. Then
    \[ \Hilb((\rS(V^*) \otimes \wedge M^*)^\chi; q, z)
        = q^{\deg \Delta_\chi} \frac{\prod_{i=1}^r (1 + zq^{-e_i^{M^*})}}
           {\prod_{i=1}^n (1 - q^{d_i})}. \]
\end{corollary}

\begin{remark}
  \Cref{cor:super_alternant_basis.2} gives Hilbert series for the invariants
  and alternants of the four algebras in \eqref{eq:four_algebras}.
  \Cref{tab:products} gives a summary. These formulas may also be derived from
  Shepler's results in \cite{MR2163712}
  Reiner--Shepler--Sommers \cite{RSS19} have recently given a related formula
  for $\Hilb((\rS(V^*) \otimes \wedge V^* \otimes \wedge V)^G; q, t, z)$
  when $G$ is a ``coincidental'' complex reflection group.
\end{remark}

\begin{table}[ht]
  \centering
  \begin{tabular}{c|c|c|c|c}
    & $\Hilb(\cA^G; q, z)$
    & $\Hilb(\cA^{{\det}_{V^*}}; q, z)$
    & $\Hilb(\cA^{{\det}_V}; q, z)$ \\
    \midrule
    $\rS(V^*) \otimes \wedge V^*$
      & $\prod_{i=1}^n \frac{1+q^{e_i} z}{1 - q^{d_i}}$
      & $\prod_{i=1}^n \frac{z+q^{e_i^*}}{1 - q^{d_i}}$
      & $q^{\sum_{i=1}^n e_i - e_i^*} \prod_{i=1}^n \frac{z + q^{e_i^*}}{1 - q^{d_i}}$ \\
    \midrule
    $\rS(V) \otimes \wedge V$
      & $\prod_{i=1}^n \frac{1+q^{e_i} z}{1 - q^{d_i}}$
      & $q^{\sum_{i=1}^n e_i - e_i^*} \prod_{i=1}^n \frac{z + q^{e_i^*}}{1 - q^{d_i}}$
      & $\prod_{i=1}^n \frac{z+q^{e_i^*}}{1 - q^{d_i}}$ \\
    \midrule
    $\rS(V) \otimes \wedge V^*$
      & $\prod_{i=1}^n \frac{1+q^{e_i^*} z}{1 - q^{d_i}}$
      & $\prod_{i=1}^n \frac{z+q^{e_i}}{1 - q^{d_i}}$
      & - \\
    \midrule
    $\rS(V^*) \otimes \wedge V$
      & $\prod_{i=1}^n \frac{1+q^{e_i^*} z}{1 - q^{d_i}}$
      & -
      & $\prod_{i=1}^n \frac{z+q^{e_i}}{1 - q^{d_i}}$
  \end{tabular}
  \vspace{0.5em}
  \caption{Product formulas for invariants and alternants for the algebras $\cA$ listed in
    \eqref{eq:four_algebras} when $M=V$.}
  \label{tab:products}
\end{table}

\section{Coinvariant harmonics and semi-invariants}\label{sec:coinvariant_alternants}
We next define the $G$-harmonic polynomial differential forms, which naturally involve
operators $\delta_i^*$ similar to the operators $d_i^*$ of the preceding
section. We prove the basic properties of the $G$-harmonics using
Hermitian forms and we show that in fact the $d_i^*$'s preserve the $G$-harmonics.
Finally, we prove our main result, \Cref{thm:super_coinvariant_basis}, giving an
explicit basis for the semi-invariant differential harmonics and coinvariants.
We continue the notation of \Cref{sec:super_alternants}.

\begin{definition}
  The \textit{coinvariant ideal} of $\rS(V^*) \otimes \wedge M^*$ is
  the ideal $\cJ_M^*$ generated by homogeneous
  $G$-invariants of positive degree. The coinvariant ideal $\cJ_M$ of
  $\rS(V) \otimes \wedge M$
  is defined analogously. When $M=V$, we write
  $\cJ^* \coloneqq \cJ_V^*$ and $\cJ \coloneqq \cJ_V$ in analogy
  with the classical coinvariant ideals $\cI^*$ and $\cI$.
\end{definition}

By Solomon's theorem \cite{MR0154929}, $\cJ^*
= (f_1, \ldots, f_n, \dif f_1, \ldots, \dif f_n)$ is the ideal generated
by $\{f_1, \ldots, f_n, \dif f_1, \ldots, \dif f_n\}$ and
$\cJ = (z_1, \ldots, z_n, \dif z_1, \ldots, \dif z_n)$. More generally,
Orlik-Solomon \cite{MR575083} show that, if $J_M = \Delta_M$ up
to a non-zero scalar, we have
\begin{align*}
  \cJ_M^* &= (f_1, \ldots, f_n, \omega_1^M, \ldots, \omega_r^M), \\
  \cJ_M &= (z_1, \ldots, z_n, \tilde{\omega}_1^M, \ldots, \tilde{\omega}_r^M).
\end{align*}

\begin{definition}
  The spaces of \textit{$G$-harmonic} elements of $\rS(V^*) \otimes \wedge M^*$ and
  $\rS(V) \otimes \wedge M$ are, respectively,
  \begin{align*}
    \cH(\rS(V^*) \otimes \wedge M^*)
      &\coloneqq \{\omega \in \rS(V^*) \otimes \wedge M^*
        : \langle \xi, \omega\rangle = 0 \text{ for all }\xi \in \cJ_M\} \\
    \cH(\rS(V) \otimes \wedge M)
      &\coloneqq \{\xi \in \rS(V) \otimes \wedge M
        : \langle \xi, \omega\rangle = 0 \text{ for all }\omega \in \cJ_M^*\}.
  \end{align*}
\end{definition}

That is, the harmonics are the orthogonal complements of the coinvariant ideals
$\cJ_M$ and $\cJ_M^*$ with respect to the natural perfect pairing
$\langle -, -\rangle$.

Recall the $G$-equivariant $F$-algebra homomorphism 
\begin{align*}
  \rS(V) \otimes \wedge M &\to \End_F(\rS(V^*) \otimes \wedge M^*) \\
  s \otimes m &\mapsto \partial_s \otimes \iota_m.
\end{align*}
from \Cref{lem:actions}.

\begin{definition}\label{def:delta_i_star}
  For $i = 1, \ldots, r$, let
    \[ \delta_i^* \colon \rS(V^*) \otimes \wedge M^* \to \rS(V^*) \otimes \wedge M^* \]
  denote the action of the basic derivation
  $\tilde{\omega}_i^M \in (\rS(V) \otimes \wedge M)^G$.
\end{definition}

\begin{lemma}
  Suppose $\cJ_M$ is generated by $j_1, \ldots, j_p$. Then
  \begin{align}
    \label{eq:harmonics_big.1}
    \cH(\rS(V^*) \otimes \wedge M^*)
      &= \{\omega \in \rS(V^*) \otimes \wedge M^* : j_i \cdot f = 0, i \in [p]\} \\
    \label{eq:harmonics_big.2}
    \rS(V^*) \otimes \wedge M^*
      &= \cH(\rS(V^*) \otimes \wedge M^*) \oplus \cJ_M^* \\
    \label{eq:harmonics_big.3}
    \cH(\rS(V^*) \otimes \wedge M^*)
      &\cong (\rS(V^*) \otimes \wedge M^*)/\cJ_M^* \text{ as bigraded $G$-modules} \\
    \label{eq:harmonics_big.4}
    \rS(V^*) \otimes \wedge M^*
      &= (\rS(V^*) \otimes \wedge M^*)^G \cH(\rS(V^*) \otimes \wedge M^*),
  \end{align}
  and likewise with $\cH(\rS(V) \otimes \wedge M)$.
  
  \begin{proof}
    When $F$ is a subfield of $\bC$ closed under complex conjugation,
    $\tau((\rS(V) \otimes \wedge M)^G) = (\rS(V^*) \otimes \wedge M^*)^G$
    and $\tau(\cJ_M) = \cJ_M^*$. The first three thus follow exactly as in the proof of
    \Cref{lem:harmonics_props.1}. Note that \eqref{eq:harmonics_big.4} merely
    asserts that the multiplication map
      \[  (\rS(V^*) \otimes \wedge M^*)^G \times \cH(\rS(V^*) \otimes \wedge M^*)
          \to \rS(V^*) \otimes \wedge M^* \]
    is surjective; it is not in general injective. Surjectivity may be proven by induction
    using \eqref{eq:harmonics_big.2}.
  \end{proof}
\end{lemma}

\begin{lemma}\label{lem:ddelta_comm}
  We have
    \[ d_i^* \delta_j^* + \delta_j^* d_i^* = \partial_{L_{i,j}} \]
  with $L_{ij} \in \rS(V)^G$.
  
  \begin{proof}
    If $F$ is a subfield of $\bC$ closed under complex conjugation, we may take
    \begin{align*}
      \tilde{\omega}_i^M
        &= \sum_j \tilde{J}_{ij}^M \otimes m_j \in (\rS(V) \otimes M)^G, \\
      \delta_i^*
        &= \sum_j \partial_{\tilde{J}_{ij}^M} \otimes \iota_{m_j}, \\
      \tilde{\omega}_i^{M^*}
        &= \sum_j \tilde{J}_{ij}^{M^*} \otimes \tau(m_j) \in (\rS(V) \otimes M^*)^G, \\
      d_i^*
        &= \sum_j \partial_{\tilde{J}_{ij}^{M^*}} \otimes \epsilon_{\tau(m_j)}.
    \end{align*}
    We may transfer the Hermitian form on $M^*$ to a Hermitian form on $M$
    defined by $(m_1, m_2) \coloneqq (\tau(m_2), \tau(m_1))$.
    Note that $\tau(m_k)(m_\ell) = \langle m_\ell, \tau(m_k)\rangle
    = (\tau(m_\ell), \tau(m_k)) = (m_k, m_\ell)$. Using \eqref{eq:iota_eps}, we
    calculate
    \begin{align*}
      d_i^* \delta_j^* + \delta_j^* d_i^*
        &= \sum_{k, \ell}
              (\partial_{\tilde{J}_{ik}^{M^*}} \partial_{\tilde{J}_{j\ell}^M})
              (\epsilon_{\tau(m_k)} \iota_{m_\ell} + \iota_{m_\ell} \epsilon_{\tau(m_k)}) \\
        &= \sum_{k, \ell} (\partial_{\tilde{J}_{ik}^{M^*} \tilde{J}_{j\ell}^M})
              (m_k, m_\ell) = \partial_{L_{ij}}
    \end{align*}
    where
    \begin{equation}\label{eq:L_ij}
      L_{ij} \coloneqq \sum_{k, \ell} (m_k, m_\ell)
       \tilde{J}_{ik}^{M^*} \tilde{J}_{j\ell}^M.
    \end{equation}
    Since $d_i^*\delta_j^* + \delta_j^*d_i^*$ is independent of the choice of
    basis $m_1, \ldots, m_r$, the same is true of $L_{ij}$. We have
      \[ \tilde{\omega}_i^M = g\tilde{\omega}_i^M = \sum_j g\tilde{J}_{ij}^M \otimes gm_j \]
    and similarly for $\tilde{\omega}_i^{M^*}$, so
    \begin{align*}
      L_{ij}
        &= \sum_{k,\ell} (gm_k, gm_\ell) g\tilde{J}_{ik}^{M^*} g\tilde{J}_{j\ell}^M \\
        &= g \sum_{k,\ell} (m_k, m_\ell) \tilde{J}_{ik}^{M^*} \tilde{J}_{j\ell}^M
          = gL_{ij}.
    \end{align*}
    The general case follows by extending scalars as in the proof of
    \Cref{lem:harmonics_props.1}.
  \end{proof}
\end{lemma}

\begin{corollary}\label{cor:di_closure}
  If $M^G = 0$, then $L_{ij}$ (as in \Cref{lem:ddelta_comm}) is homogeneous of
  positive degree and
    \[ d_i^* \in \End_F(\cH(\rS(V^*) \otimes \wedge M^*)). \]
  
  \begin{proof}
    Since $M^G = 0$, we have $(F \otimes M)^G = 0$, so $\deg_{\rS(V^*)}
    \tilde{\omega}_i^M \geq 1$. By \eqref{eq:L_ij}, $\deg L_{ij} \geq 2$,
    so $L_{ij} \in \cI$. Suppose
    $\omega \in \cH(\rS(V^*) \otimes M^*)$, or equivalently $\delta_j^*\omega = 0$
    and $\partial_{z_j}\omega = 0$. We must show
    $d_i^*\omega \in \cH(\rS(V^*) \otimes M^*)$. Since
    $\partial_{z_j}$ commutes with $d_i^*$, we have $\partial_{z_j}(d_i^*\omega) = 0$.
    By \Cref{lem:ddelta_comm}, we have
      \[ \delta_j^*(d_i^*\omega) = -d_i^*(\delta_j^*\omega) + \partial_{L_{ij}} \omega
          = 0 + 0 = 0 \]
    since $L_{ij} \in \cI$.
  \end{proof}
\end{corollary}

We now restate and prove our main result from the introduction.

\begin{theorem}\label{thm:super_coinvariant_basis}
  Suppose $J_{M^*} \mid \Delta_\chi$ and $M^G = 0$. Then the $2^r$ elements
    \[ \{d_{i_1}^* \cdots d_{i_k}^* \Delta_\chi \mid 1 \leq i_1 < \cdots < i_k \leq r\} \]
  form a basis of $\cH(\rS(V^*) \otimes \wedge M^*)^\chi$,
  and their images descend to a basis of
  $(\rS(V^*) \otimes \wedge M^*/\cJ_M^*)^\chi$.
    
  \begin{proof}
    Since $\cI \subset \cJ_M$,
    \begin{equation}\label{eq:harmonic_containment.1}
      \cH(V^*) \otimes \wedge M^* \supset \cH(\rS(V^*) \otimes \wedge M^*).
    \end{equation}
    We have $\Delta_{M^*} \in \cH(\rS(V^*) \otimes \wedge M^*)$ since
    $\delta_i^* \Delta_{M^*} = 0$ trivially.
    By \Cref{cor:di_closure}, the proposed elements belong to
    $\cH(\rS(V^*) \otimes \wedge M^*)$. The
    result hence follows from \Cref{thm:super_alternant_basis}. We may pass to the
    quotient by \eqref{eq:harmonics_big.2}.
  \end{proof}
\end{theorem}

\begin{corollary}\label{cor:super_coinvariant_basis}
  Suppose $J_{M^*} \mid \Delta_\chi$ and $M^G = 0$. Then
    \[ \Hilb((\rS(V^*) \otimes \wedge M^*/\cJ_M^*)^\chi; q, z)
        = q^{\deg \Delta_\chi}\prod_{i=1}^r (1 + zq^{-e_i^{M^*}}). \]
\end{corollary}

\begin{remark}
  By \eqref{eq:harmonic_containment.1}, \Cref{thm:super_coinvariant_basis} gives
  the surprising result that
  \begin{equation}\label{eq:harmonic_containment}
    (\cH(V^*) \otimes \wedge M^*)^\chi
    = \cH(\rS(V^*) \otimes \wedge M^*))^\chi.
  \end{equation}
\end{remark}

\section{The condition $J_M = \Delta_M$}\label{sec:Jm_Dm}
We retain the notation of the preceding sections and discuss the
condition $J_M = \Delta_M$ appearing in some of our results. In type $A$, with
precisely one exception, all irreducible examples are either one-dimensional or the defining
representation.

\begin{proposition}
  Suppose $G = \fS_n$ consists of $n \times n$ permutation matrices and
  let $M$ be an irreducible $\fS_n$-module. Then $J_M = \Delta_M$ up to a non-zero scalar
  if and only if $M$ is the
  trivial representation, the sign representation, the standard representation, or,
  when $n=4$, the unique degree $2$ representation.
  
  \begin{proof}[Proof (Sketch).]
    The irreducible $\fS_n$-representations are indexed by integer partitions
    $\lambda \vdash n$. The number of exponents of the irreducible indexed by $\lambda$
    is the degree $f^\lambda$. Since there are sub-exponentially many $\lambda \vdash n$
    and $\sum_{\lambda \vdash n} (f^\lambda)^2 = n!$, one heuristically expects
    $f^\lambda$ to grow super-exponentially in $n$, making $J_M = \Delta_M$
    quite rare. One may for instance apply the arguments in \cite[p.~15]{MR3857157}
    to make this intuition rigorous, though we omit the details.
  \end{proof}
\end{proposition}

If $G$ is a dihedral group, then every irreducible representation $M$ takes reflections to
reflections or the identity, so by \Cref{cor:Gutkin}(\ref{cor:Gutkin.pseudos}),
$J_M = \Delta_M$ up to a non-zero scalar. We now give somewhat less trivial examples.

\begin{example}\label{ex:Bn_JV}
  Let $\fB_n \subset O(n)$ be the Weyl group of type $B_n$ realized as $n \times n$
  signed permutation matrices whose non-zero entries are $\pm 1$.
  Let $\fZ_n \subset \fB_n$ be the group of diagonal matrices
  with diagonal entries $\pm 1$. Note that $\fS_n = \fB_n/\fZ_n$. Thus we may consider
  the standard representation $M$ of $\fS_n$ as an irreducible $\fB_n$-representation.
  The exponents of $M$ are the degrees in which $M$ appears in the coinvariant algebra
  $\bR[\fB_n]/\cI$, or equivalently in the harmonics $\cH(\fB_n)$. Since $M$ is
  $\fZ_n$-invariant, these occur in $\cH(\fB_n)^{\fZ_n}$. By Chevalley's result applied to
  $\fB_n$ (see \Cref{rem:chev}), multiplication gives an isomorphism
    \[ \bR[x_1, \ldots, x_n]^{\fB_n} \otimes \cH(\fB_n)
        \stackrel{\sim}{\to} \bR[x_1, \ldots, x_n]. \]
  Taking $\fZ_n$-invariants and simplifying,
    \[ \bR[x_1^2, \ldots, x_n^2]^{\fS_n} \otimes \cH(\fB_n)^{Z_n}
        \stackrel{\sim}{\to} \bR[x_1^2, \ldots, x_n^2]. \]
  Again applying Chevalley's result to $\fS_n$, it follows that $\cH(\fB_n)^{\fZ_n}
  \cong \cH(\fS_n)$ where $x_i^2 \mapsto x_i$. In particular, $J_M = \Delta_M$
  up to a non-zero multiple, and both have degree $2 + 4 + \cdots + 2(n-1) =
  2\binom{n}{2} = n(n-1)$.
\end{example}

\begin{example}
  Let $\fD_n \subset \fB_n \subset O(n)$ be the Weyl group of type $\fD_n$ consisting of
  signed permutation matrices with evenly many negative entries. Let $\fZ_n'
  = \fD_n \cap \fZ_n$. Note that $\fB_n/\fZ_n' \cong \fS_n$. It is easily checked that
    \[ \bR[x_1, \ldots, x_n]^{\fZ_n'} = \bR[x_1^2, \ldots, x_n^2]
        \oplus \bR[x_1^2, \ldots, x_n^2] x_1 \cdots x_n \]
  and
    \[ \bR[x_1, \ldots, x_n]^{\fD_n} = \bR[x_1^2, \ldots, x_n^2]^{\fS_n}
        \oplus \bR[x_1^2, \ldots, x_n^2]^{\fS_n} x_1 \cdots x_n. \]
  Using the same argument as in \Cref{ex:Bn_JV}, the standard representation $M$ of
  $\fS_n$ yields an irreducible representation of $\fD_n$ with $J_M = \Delta_M$
  up to a non-zero multiple, and both have degree $n(n-1)$.
\end{example}

The Weyl group $\fF_4$ of type $F_4$ is isomorphic with $\fS_3 \ltimes \fD_4$.
One may then view the standard representation $M$ of $\fS_3$
as an $\fF_4$-module and check that the exponents of $M$ are $4, 8$
and $J_M = \Delta_M$ up to a non-zero scalar. The hypotheses of
\Cref{cor:Gutkin}(\ref{cor:Gutkin.pseudos}) are satisfied in this case.
Further examples exist. For example, Orlik--Solomon observe that Galois
conjugates of $V$ and $V^*$ have this property. A complete classification is
not known.

\section{Acknowledgements}
The first named author would like to thank John Mahacek, Vic Reiner, Brendon Rhoades,
Bruce Sagan, Anne Shepler, and Mike Zabrocki for helpful discussions on this and related
projects, and for generously sharing their preprints. The second named author thanks Adriano
Garsia for pointing out that our result for the signum representation of
$S_n$ can be derived from the Zabrocki conjecture.

\bibliography{refs}{}
\bibliographystyle{acm}

\end{document}